\documentclass[11pt]{article}
\usepackage{amssymb,amsmath,amsfonts,amsthm,array,bm,bbm,color,graphicx,hyperref}

\hypersetup{colorlinks=true, linkcolor=red, anchorcolor=blue, citecolor=blue}
\providecommand{\U}[1]{\protect\rule{.1in}{.1in}}

\numberwithin{equation}{section}

\newtheorem{theorem}{Theorem}[section]
\newtheorem{lemma}[theorem]{Lemma}
\newtheorem{corollary}[theorem]{Corollary}

\newtheorem{definition}[theorem]{Definition}

\def\<{\langle}
\def\>{\rangle}
\def\d{{\rm d}}

\def\R{\mathbb{R}}
\def\T{\mathbb{T}}

\def\p{\partial}

\begin{document}

\title{Mean-field derivation of Landau-like equations}


\author{José Antonio Carrillo\footnote{Email: carrillo@maths.ox.ac.uk. Mathematical Institute, University of Oxford, Oxford, OX2 6GG, UK.} \,   Shuchen Guo\footnote{Email: guo@maths.ox.ac.uk. Mathematical Institute, University of Oxford, Oxford, OX2 6GG, UK.} \,   Pierre-Emmanuel Jabin\footnote{Email: pejabin@psu.edu. Department of Mathematics and Huck Institutes, Pennsylvania State University, State College, PA 16801, USA.}}

\maketitle

\vspace{-20pt}

\begin{abstract}
We derive  a class of  space homogeneous Landau-like equations from stochastic interacting particles. Through the use of relative entropy,
we obtain quantitative bounds on the distance between the solution of the N-particle Liouville equation and the
tensorised solution of the limiting Landau-like equation.
\end{abstract}







\section{Introduction}
We consider Landau-like equations on the torus $\T^d$ that
\begin{equation}\label{Landau-like}
\p_t f  =\nabla \cdot \big[(a\ast f)\nabla f-(b\ast f) f\big],
\end{equation}
where the matrix-valued function $a$ is symmetric and uniformly bounded from above and below in the sense of bilinear form as
$$
\lambda_1\operatorname{Id}\leq a \leq \lambda_2\operatorname{Id}, \quad 0<\lambda_1\leq \lambda_2.
$$
Here $\operatorname{Id}$ is the identity matrix;
and the vector field $b=\nabla\cdot a$ as well as its divergence $\nabla\cdot b$ is bounded. We also assume the initial data of \eqref{Landau-like} satisfies $f^0\in W^{2,\infty}(\T^d)$. 

Equation \eqref{Landau-like} can be written in non-divergence form as
$$
\p_t f=(a\ast f):\nabla^2 f-\nabla\cdot(b\ast f)f.
$$
From the classical theory of advection-diffusion equations, we can assume the solution of \eqref{Landau-like} belongs to $f\in L^{\infty}([0,T], W^{2,\infty}(\T^d))$ with $\int_{\T^d}f=1$ and $\inf_{v\in \T^d} f(v)>0$ for all $t\in[0,T]$.

The Landau equation \cite{L1936,V02} plays an important role in kinetic theory, and in particular to model a plasma of charged particles. It can be formally derived from Boltzmann equation, in which grazing collision prevails.  The true space homogeneous Landau equation on $\R^d$ has the same structure as \eqref{Landau-like} but with matrix-valued function $a$ and vector-valued function $b$ as
$$
a(z)=|z|^{\gamma}\Big(\operatorname{Id}|z|^{2}-z\otimes z\Big) \quad\text{and}\quad b(z)=\nabla\cdot a(z)=-2|z|^\gamma z,
$$
where $\gamma\geq-d$. 
The solution $f(t,v)$ corresponds to the probability of finding a particle in the plasma at time $t$ with velocity $v$.  In the case $d=3$, one usually speaks of hard potentials when $\gamma\in(0,1]$, 
Maxwellian
potential when $\gamma=0$, moderately soft potentials when $\gamma \in[-2,0)$, very soft potential when $\gamma \in[-3,-2)$,
and the special case of Coulomb potential  corresponding to $\gamma=-3$. Our assumptions for Landau-like equation on $a$ and $b$ avoid the possible degeneracy and singularity at the origin, but keep the structure $b=\nabla\cdot a$. 

In terms of the properties of the Landau equation,  for hard potentials, well-posedness, regularity and
large-time behavior have been studied by Desvillettes–Villani \cite{DV00i,DV00ii} and Fournier-Heydecker \cite{FH21}; for Maxwellian case, these are given by Villani \cite{V98i}; for  moderately soft potentials, a global well-posedness result is obtained in \cite{FG09}; for very soft potentials, \cite{V98ii} defines the H-solution and proves its existence, but the regularity and uniqueness of H-solutions remain open. Recently, Guillen and Silvestre showed that the classical solution of Landau equation will not blow up for all $\gamma\in[-3,1]$ \cite{GS23}.

The rigorous derivation of the Landau equation directly from many-particle systems is still a challenge. Using Newtonian dynamics, Kac proposed to derive the Boltzmann equation from stochastic particle system in sense of mean field limit, and he gave the mathematical definition of  molecular chaos \cite{K56}. For a detailed discussion about mean field limit and propagation of (molecular) chaos, one can see \cite{S91ii,JW17}. We will adopt this idea to derive space homogeneous Landau-like equations from coupled systems of SDEs for the particles, while proving the propagation of chaos.

In the case of Maxwellian potential,  Fontbona, Gu\'{e}rin and M\'{e}l\'{e}ard \cite{FGM09} obtained a quantitative propagation of chaos for Landau-like equations, and  Fournier \cite{F09} improved the rate of convergence. Later,  Carrapatoso was able to prove a uniform in time quantitative propagation of chaos in \cite{C16}. When $\gamma\in[0,1]$, Fournier and Guillin also derived a quantitative result \cite{FG17}. For soft potential, Fournier and Hauray deal with both  $\gamma\in(-1,0)$ and $\gamma\in(-2,-1]$ in \cite{FH16}, for the former case, they obtain a  rate of convergence; for the latter case, they also prove the propagation of chaos but without an explicit rate. All these results are proven by using coupling techniques.

Inspired by the work of Jabin and Wang \cite{JW16,JW18}, we prove a quantitative propagation of chaos by controlling the relative entropy, which yields the derivation of Landau-like equations from stochastic particle systems. This work is organised as follows. Section \ref{main result} is dedicated to introduce our particle systems and state main results; and the proof of Theorem \ref{poc} is given in Section \ref{proof}.
\section{Main result}\label{main result}
We consider the following stochastic $N$-particle systems on $\T^d$: 
$$
\d V_t^i=\frac{2}{N}\sum_{j=1}^{N}b(V_t^i-V_t^j)\d t+\sqrt{2}\Big(\frac{1}{N}\sum_{j=1}^{N}a(V_t^i-V_t^j)\Big)^{\frac{1}{2}}\d B_t^i,
$$
where $(B_t^i)_{i\geq 1}$ are i.i.d.  $d$-dimensional Brownian motions and the diffusion coefficient matrix is
a unique square root of the nonnegative symmetric matrix. We use the convention that $a(0)=0$ and $b(0)=0$ to  omit the notation $i\neq j$. We notice that, under our assumptions on the particle system, the particles are exchangeable, thus we assume that the initial  joint distribution of $(V_0^1,\ldots,V_0^N)$ is a symmetric probability measure $f_N(0)$. 

The existence and uniqueness of strong solution to the particle systems (SDEs) have been proved in \cite{F09}.   Applying It\^{o}'s formula and the relation $\nabla\cdot a=b$, we can derive the evolution (Liouville equation) of $N$-particles joint distribution $f_N(t, V)$, $V=(v^1,\ldots,v^N)$ on $\T^{dN}$ as 
\begin{equation}\label{fn}
\begin{aligned}
\p_t f_N= &\, \sum_{i=1}^{N}\nabla_{v^{i}}\cdot\Big[\frac{1}{N}\sum_{j=1}^{N}a(v^{i}-v^{j})\nabla_{v^{i}}f_N-\frac{1}{N}\sum_{j=1}^{N}b(v^{i}-v^{j})f_N\Big],
\end{aligned}
\end{equation}
where the 
initial value is $f_N(0)$. There exists at least one entropy solution  of  \eqref{fn} defined as follows \cite{JW18}.
\begin{definition}[Entropy solution]
For $t\in [0,T]$, 
  a density function $f_N\in\T^{dN}$, with $f_N\geq 0$ and $\int_{\T^{dN}}f_N=1$, is called an entropy solution of  \eqref{fn}  if and only if 
$$
\begin{aligned}
 \int_{\T^d}f_N(t)\log f_N(t) \d V&+\sum_{i=1}^N\int_0^t\int_{\T^d}f_N\frac{1}{N}\sum_{j=1}^{N}a(v^{i}-v^{j}):\nabla_{v^{i}}\log f_N\otimes\nabla_{v^{i}}\log f_N\d V\\
&+\sum_{i=1}^N\int_0^t\int_{\T^d}f_N\frac{1}{N}\sum_{j=1}^{N}\nabla\cdot b(v^{i}-v^{j})\d V\leq  \int_{\T^d}f_N(0)\log f_N(0) \d V.
\end{aligned}
$$
\end{definition}
To prove the propagation of chaos, we aim to estimate the distance between the solution of Liouville equation $f_N$ and tensorised solution of the Landau equation $\bar f_N=f^{\otimes N}$ in the sense of relative entropy,
\begin{equation*}
H_N(f_N|\bar f_N)=\frac{1}{N}\int_{\T^{dN}}f_N\log\frac{f_N}{\bar f_N}\d V=\frac{1}{N}\int_{\T^{dN}}f_N\log f_N\d V-\frac{1}{N}\int_{\T^{dN}}f_N\log \bar f_N\d V.
\end{equation*}
Its subadditivity implies the bound on the distance between $k$-th marginals
\begin{equation*}
f_{k,N}(t,v^1,\ldots,v^k)=\int_{\T^{d(N-k)}}f_N(t,v^1,\ldots,v^N)\d v^{k+1}\ldots \d v^N,
\end{equation*}
and tensorised
$f^{\otimes k}$ as
\begin{equation}\label{subaddivity}
H_k(f_{k,N}|f^{\otimes k})=\frac{1}{k}\int_{\T^{dk}}f_{k,N}\log\frac{f_{k,N}}{f^{\otimes k}}\d v^1\cdots\d v^k\leq H_N(f_N|f^{\otimes N}).
\end{equation}
For simplicity, we denote the $H_N(f_N(t)|\bar f_N(t))$ as $H_N(t)$. Now we state our main result.
\begin{theorem}\label{poc}
Under assumptions above, 
there exists some positive constant $C_1$  and $C_2$ independent with $N$ such that the relative entropy of $f_N$ and $\bar f_N$ on the torus $\T^{dN}$ has the following estimate
$$
H_N(t)\leq \Big(H_N(0)+\frac{C_1}{N}\Big)e^{C_2t}.
$$
\end{theorem}

Then, Theorem \ref{poc} implies the quantitative propagation of chaos result:
\begin{corollary}\label{corollary}
    Under assumptions above, and further assuming $\sup_NNH_N(0)<\infty$, one has the strong propagation of chaos, for some   constant $C_3$ independent with $N$,
    $$
    \|f_{k,N}-f^{\otimes k}\|_{L^{\infty}([0,T],L^1(\T^{dk}))}\leq \frac{C_3}{\sqrt{N}}.
    $$
\end{corollary}
The proof of Corollary \ref{corollary} is straightforward by applying \eqref{subaddivity} and Csisz\'{a}r-Kullback-Pinsker inequality as for any functions $g_1$ and $g_2$ on $\T^{dk}$ as
$$
\|g_1-g_2\|_{L^1(\T^{dk})}\leq \sqrt{2kH_k(g_1|g_2)}.
$$
\section{Proof of Theorem \ref{poc}}\label{proof}
We firstly derive the evolution of relative entropy $H_N$, and it holds
$$
\begin{aligned}
\frac{\d}{\d t}\left(\frac{1}{N}\int_{\T^{dN}}f_N\log f_N\d V\right)=&\,\frac{1}{N}\int_{\T^{dN}}(1+\log f_N)\p_t f_N\d V
\\=&-\frac{1}{N}\sum_{i=1}^{N}\int_{\T^{dN}}\frac{1}{N}\sum_{j=1}^{N}a(v^{i}-v^{j}):\frac{\nabla_{v^{i}}f_N}{\sqrt{f_N}}\otimes\frac{\nabla_{v^{i}}f_N}{\sqrt{f_N}}\d V
\\&+\frac{1}{N}\sum_{i=1}^{N}\int_{\T^{dN}}\frac{1}{N}\sum_{j=1}^{N}b(v^{i}-v^{j})\cdot\nabla_{v^{i}}f_N\d V,
\end{aligned}
$$
where we plug in \eqref{fn} and integrate by parts in the second step; similarly, we have
$$
\begin{aligned}
\frac{\d}{\d t}\left(\frac{1}{N}\right.&\left.\!\!\int_{\T^{dN}}f_N\log \bar f_N\d V\right)
=\frac{1}{N}\int_{\T^{dN}}\Big(\log \bar f_N \p_t f_N + f_N \frac{ \p_t \bar f_N}{\bar f_N}\Big)\d V
\\=&\,-\frac{1}{N}\sum_{i=1}^{N}\int_{\T^{dN}}\!\!\frac{1}{N}\sum_{j=1}^{N}a(v^{i}-v^{j}):
\frac{\nabla_{v^{i}}\bar f_N\otimes\nabla_{v^{i}}f_N}{\bar f_N} \d V\\&\,+\frac{1}{N}\sum_{i=1}^{N}\int_{\T^{dN}}\!\!\frac{1}{N}\sum_{j=1}^{N}b(v^{i}-v^{j})\cdot\frac{\nabla_{v^{i}}\bar f_N}{\bar f_N}f_N\d V
\\&\,-\frac{1}{N}\sum_{i=1}^{N}\int_{\T^{dN}} a\ast f(v^{i}) :\nabla_{v^{i}}\frac{f_N}{\bar f_N}\otimes\nabla_{v^{i}}\bar{f}_N\d V+\frac{1}{N}\sum_{i=1}^{N}\int_{\T^{dN}}b\ast f(v^{i})\cdot \nabla_{v^{i}}
\frac{ f_N }{\bar f_N} \bar f_N\d V.
\end{aligned}
$$
Using the identity 
$$
\nabla_{v^i} \frac{f_N}{\bar{f}_N}=\frac{\bar{f}_N \nabla_{v^i} f_N-f_N \nabla_{v^i} \bar{f}_N}{\bar{f}_N^2},
$$
enables us to rewrite 
\begin{equation}\label{RE}
\begin{aligned}
\frac{\d}{\d t}H_N(t)= &\,\frac{\d}{\d t}\left(\frac{1}{N}\int_{\T^{dN}}f_N\log f_N\d V\right)-\frac{\d}{\d t}\left(\frac{1}{N}\int_{\T^{dN}}f_N\log \bar f_N\d V\right)
\\=&\, -\frac{1}{N}\sum_{i=1}^{N}\int_{\T^{dN}}f_N\frac{1}{N}\sum_{j=1}^{N}a(v^{i}-v^{j}):\nabla_{v^{i}}\log\frac{f_N}{\bar f_N}\otimes\nabla_{v^{i}}\log\frac{f_N}{\bar f_N}\d V\\
 &\, +\frac{1}{N}\sum_{i=1}^{N}\int_{\T^{dN}} f_N\Big[a\ast f(v^{i})-\frac{1}{N}\sum_{j=1}^{N}a(v^{i}-v^{j})
\Big]:\nabla_{v^{i}}\log\frac{f_N}{\bar f_N}\otimes\nabla_{v^{i}}\log\bar f_N\d V
\\ &\, -\frac{1}{N}\sum_{i=1}^{N}\int_{\T^{dN}}f_N\Big[b\ast f(v^{i})-\frac{1}{N}\sum_{j=1}^{N}b(v^{i}-v^{j})\Big]\cdot \nabla_{v^{i}}\log
\frac{ f_N }{\bar f_N} \d V\\
=:&\, I_1+I_2+I_3.
\end{aligned}
\end{equation}

The first term on the right-hand side of \eqref{RE} can be bounded by the assumption on the minimal eigenvalue of matrix $a$ as
\begin{equation*}
  \begin{aligned}
I_1=& \,-\int_{\T^{dN}}f_N\frac{1}{N}\sum_{j=1}^{N} a(v^{1}-v^{j}):\nabla_{v^{1}}\log\frac{f_N}{\bar f_N}\otimes\nabla_{v^{1}}\log\frac{f_N}{\bar f_N}\d V\\
\leq &\,-\int_{\T^{dN}}f_N\,\lambda_1\operatorname{Id}:\nabla_{v^{1}}\log\frac{f_N}{\bar f_N}\otimes\nabla_{v^{1}}\log\frac{f_N}{\bar f_N}\d V=-\lambda_1\int_{\T^{dN}}f_N\Big|\nabla_{v^{1}}\log\frac{f_N}{\bar f_N}\Big|^2\d V,
\end{aligned}  
\end{equation*}
where we used $v^1$ instead of averaging all $v^i$ by exchangeability. Simple inequality $xy\leq \frac{1}{\lambda_1}x^2+\frac{\lambda_1}{4} y^2$ for positive constant $\lambda_1$ implies the estimate on the last two  terms in \eqref{RE} as
$$
\begin{aligned}
I_2+I_3=&\int_{\T^{dN}} f_N  \Big[a\ast f(v^{1})-\frac{1}{N}\sum_{j=1}^{N}a(v^{1}-v^{j})
\Big]:\nabla_{v^{1}}\log\frac{f_N}{\bar f_N}\otimes\nabla_{v^{1}}\log\bar f_N\d V
\\ 
&\,-\int_{\T^{dN}}f_N \Big[b\ast f(v^{1})-\frac{1}{N}\sum_{j=1}^{N}b(v^{1}-v^{j})\Big]\cdot \nabla_{v^{1}}\log
\frac{ f_N }{\bar f_N} \d V  \\
\leq &\, \frac{\lambda_1}{ 2}\int_{\T^{dN}}f_N\Big|\nabla_{v^{1}}\log\frac{f_N}{\bar f_N}\Big|^2\d V\\&\,+\frac{1}{\lambda_1}\int_{\T^{dN}}f_N\Big\|a\ast f(v^{1})-\frac{1}{N}\sum_{j=1}^{N}a(v^{1}-v^{j})\Big\|^2\bigg|\frac{\nabla_{v^{1}}\bar f_N}{\bar f_N}\bigg|^2\d V\\&\,+\frac{1}{\lambda_1}\int_{\T^{dN}}f_N\Big|b\ast f(v^{1})-\frac{1}{N}\sum_{j=1}^{N}b(v^{1}-v^{j})\Big|^2\d V,
\end{aligned}
$$
where we take the Frobenius norm for matrices. Also we notice that under our assumption on the solutions $f(t,v)$ of the Landau-like equation \eqref{Landau-like}, it holds
$$
\sup_{V\in \T^{dN},t\in[0,T]}\bigg|\frac{\nabla_{v^{1}}\bar f_N}{\bar f_N}\bigg|=\sup_{v\in \T^{d},t\in[0,T]}\bigg|\frac{\nabla_{v}f(t,v)}{ f(t,v)}\bigg|<\infty.
$$
Combining the estimates above, we get the bound of \eqref{RE} that 
\begin{equation*}
\begin{aligned}
\frac{\d}{\d t}H_N(t)\leq &\, \frac{C_f}{\lambda_1}\sum_{\alpha,\beta=1}^d\int_{\T^{dN}}f_N\Big(a_{\alpha,\beta}\ast f(v^{1})-\frac{1}{N}\sum_{j=1}^{N}a_{\alpha,\beta}(v^{1}-v^{j})\Big)^2\d V\\&+\frac{1}{\lambda_1}\sum_{\alpha=1}^d\int_{\T^{dN}}f_N\Big(b_{\alpha}\ast f(v^{1})-\frac{1}{N}\sum_{j=1}^{N}b_{\alpha}(v^{1}-v^{j})\Big)^2\d V.
\end{aligned}
\end{equation*}
The following lemma is the same as  \cite[Lemma 1]{JW18}.
\begin{lemma}
 For any two probability densities $f_N$ and $\bar{f}_N$ on $\T^{dN}$, and any $\Phi \in L^{\infty}\left(\T^{dN}\right)$, one has that $\forall \eta>0$,
$$
\int_{\T^{dN}}f_N \Phi  \d V \leq \frac{1}{\eta}\left(H_N\left(f_N \mid \bar f_N\right)+\frac{1}{N} \log \int_{\T^{dN}} \bar{f}_N e^{N \eta \Phi} \d V\right) .
$$
\end{lemma}
Applying this lemma with bounded functions$$
\Phi_1=\Big(a_{\alpha,\beta}\ast f(v^{1})-\frac{1}{N}\sum_{j=1}^{N}a_{\alpha,\beta}(v^{1}-v^{j})\Big)^2, \quad
\Phi_2=\Big(b_{\alpha}\ast f(v^{1})-\frac{1}{N}\sum_{j=1}^{N}b_{\alpha}(v^{1}-v^{j})\Big)^2
$$
respectively, we deduce that
$$
\begin{aligned}
\frac{\d}{\d t}H_N(t)\leq &\,\frac{C_fd^2+d}{\lambda_1\eta} H(t)
\\&\,+\frac{C_f}{\lambda_1\eta N}\sum_{\alpha,\beta=1}^d\log\int_{\T^{dN}}\bar f_N \exp\bigg\{\eta N\Big(a_{\alpha,\beta}\ast f(v^{1})-\frac{1}{N}\sum_{j=1}^{N}a_{\alpha,\beta}(v^{1}-v^{j})\Big)^2\bigg\}\d V\\&\,+\frac{1}{\lambda_1\eta N}\sum_{\alpha=1}^d\log\int_{\T^{dN}}\bar f_N \exp\bigg\{\eta N\Big(b_{\alpha}\ast f(v^{1})-\frac{1}{N}\sum_{j=1}^{N}b_{\alpha}(v^{1}-v^{j})\Big)^2\bigg\}\d V.
\end{aligned}
$$
And we notice that the identity holds
$$
\eta N\Big(a_{\alpha,\beta}\ast f(v^{1})-\frac{1}{N}\sum_{j=1}^{N}a_{\alpha,\beta}(v^{1}-v^{j})\Big)^2=\Big(\frac{1}{\sqrt{N}}\sum_{j=1}^{N}\sqrt{\eta}\big(a_{\alpha,\beta}\ast f(v^{1})-a_{\alpha,\beta}(v^{1}-v^{j})\big)\Big)^2;
$$
similarly, it has
$$
\eta N\Big(b_{\alpha}\ast f(v^{1})-\frac{1}{N}\sum_{j=1}^{N}b_{\alpha}(v^{1}-v^{j})\Big)^2=\Big(\frac{1}{\sqrt{N}}\sum_{j=1}^{N}\sqrt{\eta}\big(b_{\alpha}\ast f(v^{1})-b_{\alpha}(v^{1}-v^{j})\big)\Big)^2.
$$
To give further estimate, we will take advantage of \cite[Theorem 3]{JW18} as follows.
\begin{lemma}\label{cancel} Assume that a scalar function $\psi \in L^{\infty}$ with $\|\psi\|_{L^{\infty}}<\frac{1}{2 e}$, and that for any fixed $z$, $\int_{\T^d} \psi(z, v) f(v) \d v=0$, then
$$
\begin{aligned}
&\int_{\T^{dN}} \bar{f}_N \exp \Big(\frac{1}{N} \sum_{j_1, j_2=1}^N \psi\left(v^{1}, v^{j_1}\right) \psi\left(v^{1}, v^{j_2}\right)\Big) \d V< C_0.
\end{aligned}
$$
\end{lemma}
For each entry $\alpha,\beta$, denoting functions $\psi_1$ and $\psi_2$ respectively as
$\psi_1(z,v)=\sqrt{\eta}\big(a_{\alpha,\beta}\ast f(z)-a_{\alpha,\beta}(z-v)\big)$ and $\psi_2(z,v)=\sqrt{\eta}\big(b_{\alpha}\ast f(z)-b_{\alpha}(z-v)\big)$,
we can choose some suitable $\eta$ such that each component of $\psi_1$ and $\psi_2$ satisfying the assumption that $\|\psi\|_{L^\infty}<\frac{1}{2e}$ in Lemma \ref{cancel}. Then it holds that $$
\frac{\d}{\d t}H_N(t)\leq\frac{C_fd^2+d}{\lambda_1\eta} H(t)
+\frac{C_0'}{\lambda_1\eta N}. 
$$
Therefore, Gronwall's lemma implies the main result Theorem \ref{poc}.

\section*{Acknowledgements}
JAC was supported by the Advanced Grant Nonlocal-CPD (Nonlocal PDEs for Complex Particle Dynamics: Phase Transitions, Patterns and Synchronization) of the European Research Council Executive Agency (ERC) under the European Union’s Horizon 2020 research and innovation programme (grant agreement No.~883363).

\end{document}